\documentclass{article}
\usepackage[english]{babel}
\usepackage{amsmath}
\usepackage{amsthm}
\usepackage{amssymb}
\usepackage{amscd}
\usepackage{pst-all}
\usepackage[latin1]{inputenc}

\newtheorem{teor}{Theorem}[section]
\newtheorem{defi}{Definition}
\newtheorem{lema}[teor]{Lemma}
\newtheorem{prop}[teor]{Proposition}
\newtheorem{cor}[teor]{Corollary}
\newtheorem{rem}[teor]{Remark}
\newtheorem{ejem}[teor]{Example}

\begin{document}

\title{Abelian exact subcategories closed under predecessors}

\author{Ibrahim Assem\\
Departement de Mathématique et Informatique\\ Université de
Sherbrooke\\ Sherbrooke, Québec J1K 2R1\\
CANADA\\ {\it ibrahim.assem@usherbrooke.ca} \vspace*{0.5cm} \and Manuel Saorín\\ Departamento de Matemáticas\\
Universidad de Murcia, Aptdo. 4021\\
30100 Espinardo, Murcia\\
SPAIN\\ {\it msaorinc@um.es}}

\date{Dedicated to Claus M. Ringel on the occasion of his 60th birthday}

\thanks{The first author gratefully acknowledges partial support from the NSERC of Canada.
The second author thanks the D.G.I. of the Spanish Ministry of
Science and Technology and the Fundación "Séneca" of Murcia for
their financial support}

\maketitle

\begin{abstract}

{\bf In the category of finitely generated  modules over an
artinian ring, we classify all the abelian exact subcategories
closed under predecessors or, equivalently, all the split torsion
pairs with torsion-free class closed under quotients. In the
context of Artin algebras, the result is then applied to the left
part of the module category and to local extensions of hereditary
algebras
 }

\end{abstract}

\section{Introduction}
Let $A$ be an Artin algebra, $mod_A$ be its category of right
$A$-modules,  and $ind_A$ be a full subcategory of $mod_A$
consisting of a complete set of representatives of the isomorphism
classes of indecomposable $A$-modules. The left part
$\mathcal{L}_A$ and the right part $\mathcal{R}_A$ of $mod_A$ were
introduced by Happel, Reiten and Smal$\emptyset$ in their study of
quasi-tilted algebras (\cite{HRS}). These have repeatedly proved
their usefulness in the study of homological properties of the
algebra. Our initial motivation for the present paper was the
following question: when is the additive closure
$add(\mathcal{L}_A)$ of $\mathcal{L}_A$ an abelian exact
subcategory of $mod_A$? (see definition below). As our study
advanced, we noticed that the particular consideration of
$\mathcal{L}_A$ was not essential, and our goal then shifted to
classify all the full subcategories $\mathcal{C}\subseteq ind_A$,
closed under predecessors,  such that $add(\mathcal{C})$ is an
abelian exact subcategory of $mod_A$. This is easily seen to be
equivalent to
 the classification of all split torsion pairs in $mod_A$,  with
torsion-free class closed under quotients. In addition, we
realized that the restriction to Artin algebras was not necessary
and that our classification held in the more general context of
(right) artinian rings. The desired classification is given in
corollary \ref{clasification} as a direct consequence of our main
result, theorem \ref{main-theorem}.

This theorem states that, for a basic and connected right artinian
ring $A$, the existence of such a subcategory $\mathcal{C}$ of
$ind_A$ is equivalent to the existence of an isomorphism $A\cong
\begin{pmatrix} C & 0\\ M & B \end{pmatrix}$, where
$M$ is a $B-C-$bimodule which is hereditary injective over $C$,
and such that $\mathcal{C}$ gets identified with $ind_C$. In case
$A$ is an Artin algebra or, more generally, an artinian ring with
selfduality,  our methods can be dualized to yield a
classification of those subcategories $\mathcal{C}\subseteq ind_A$
closed under successors and such that $add(\mathcal{C})$ is an
abelian exact subcategory of $mod_A$. We leave the primal-dual
translation to the reader.

The paper is organized as follows. Section 2 is devoted to proving
the main theorem, for which we
 need several equivalent characterizations of the desired
 subcategories (see
Proposition \ref{auxilio-fundamental} below). Section 3 contains
applications of the theorem to Artin algebras, in the case where
 $\mathcal{C}=\mathcal{L}_A$. Thus, we prove that if the
quiver of $A$ has no oriented cycles, then $add(\mathcal{L}_A)$ is
an abelian exact subcategory of $mod_A$ if and only if $A$ is
hereditary (see Corollary \ref{LA} below). We also prove that if
$A$ is a local extension of a hereditary algebra $H$  (by a
bimodule $_RM_H$), then $add(\mathcal{L}_A)$ is an abelian exact
subcategory of $mod_A$ if, and only if, $M_H$ is injective (see
Proposition \ref{final}).

\section{The main theorem}

Throughout this section, $A$ is a basic right artinian ring, which
we assume  connected (that is,  indecomposable as a  ring).
Modules are finitely generated right modules. All subcategories of
 $mod_A$ or $ind_A$  are assumed closed under isomorphic images.
For a full subcategory  $\mathcal{C}$  of $mod_A$, we  denote by
$add(\mathcal{C})$ the full subcategory of $mod_A$ having as
objects the direct summands of finite direct sums of modules in
$\mathcal{C}$. We also write briefly $X\in\mathcal{C}$ to express
that $X$ is an object of $\mathcal{C}$ . For an $A$-module $M$,
$Gen(M)$  stands for the full subcategory consisting of those
modules which are generated by $M$, that is, which are quotients
of modules in $add(M)$. We refer the reader to \cite{AF} and
\cite{ARS}[Chapter I] for  concepts about artinian rings not
specifically defined here.

Given  $X,Y\in ind_A$, a {\bf path} from $X$ to $Y$  is a sequence
$X=X_0\stackrel{f_1}{\longrightarrow} X_1\rightarrow
...\stackrel{f_t}{\longrightarrow}X_t=Y$ of non-zero morphisms
$f_i$ between indecomposable $A$-modules. In this case, we say
that $X$ is a {\bf predecessor} of $Y$ (and that $Y$ is a {\bf
successor} of $X$). A full subcategory $\mathcal{C}\subseteq
ind_A$ is called {\bf closed under predecessors} when every
predecessor of a module in $\mathcal{C}$ lies in $\mathcal{C}$.
When $\mathcal{C}$ is closed under predecessors, the direct sum
$P=P_\mathcal{C}$ of all (indecomposable) projective modules in
$\mathcal{C}$ is called the {\bf supporting projective} module of
$\mathcal{C}$.

We recall that a pair $(\mathcal{T},\mathcal{F})$ of full
subcategories of $mod_A$ is called a {\bf torsion pair}, when it
satifies the following two conditions: i) a module $X_A$ is in
$\mathcal{T}$ if, and only if, $Hom_A(X,F)=0$ for all
$F\in\mathcal{F}$, and ii) a module $X_A$ is in $\mathcal{F}$ if,
and only if, $Hom_A(T,X)=0$ for all $T\in\mathcal{T}$. In this
case, we have an idempotent subfunctor of the identity
$t:mod_A\longrightarrow mod_A$, called the {\bf torsion radical},
such that $X\in\mathcal{T}$ if and only if  $t(X)=X$. The class
$\mathcal{T}$ is  called the {\bf torsion class}, and the class
$\mathcal{F}$ is called the {\bf torsion-free class} of the pair.
The pair $(\mathcal{T},\mathcal{F})$ is called {\bf split} when
$t(X)$ is a direct summand of $X$, for all $X\in mod_A$, or,
equivalently, when every indecomposable $A$-module $X$ either
belongs to $\mathcal{T}$ or to $\mathcal{F}$.

 The following lemma is well-known.

\begin{lema} \label{split-torsion}
Let $\mathcal{C}$ be a full subcategory of $ind_A$. The following
assertions are equivalent:

\begin{enumerate}
\item[(1)] $\mathcal{C}$ is closed under predecessors. \item[(2)]
If $X\in ind_A$ then either $X\in\mathcal{C}$ or  $Hom_A(X,Y)=0$,
for all $Y\in\mathcal{C}$. \item[(3)] $add(\mathcal{C})$ is the
torsion-free class of a split torsion pair in $mod_A$.
\end{enumerate}
\end{lema}

In this paper we use the following terminology.

\begin{defi}
A full subcategory $\mathcal{A}$ of $mod_A$ is said to be an {\bf
abelian exact subcategory}, when it is abelian as a category and
the inclusion functor $\mathcal{A}\hookrightarrow mod_A$ is exact
\end{defi}

It is easily seen that a full subcategory $\mathcal{A}$ is an
abelian exact subcategory of $mod_A$ if, and only if, it is closed
under kernels and cokernels. In general, a full subcategory can be
abelian as a category without being an abelian exact subcategory
of $mod_A$.

\begin{prop} \label{auxilio-fundamental}
Let $\mathcal{C}$ be a full subcategory of $ind_A$ closed under
predecessors. The following statements are equivalent:

\begin{enumerate}
\item[(1)] $add(\mathcal{C})$ is an abelian exact subcategory of
$mod_A$. \item[(2)]$add(\mathcal{C})$ is closed under cokernels.
\item[(3)]$add(\mathcal{C})$ is closed under quotients. \item[(4)]
For every (indecomposable) projective $P_0\in\mathcal{C}$, we have
$top(P_0)\in\mathcal{C}$. \item[(5)]$\mathcal{C}$ is closed under
composition factors.
 \item[(6)]$add(\mathcal{C})=Gen(P)$, where $P$ is the supporting projective
module of $\mathcal{C}$.
\end{enumerate}
\end{prop}
\begin{proof}

Since $\mathcal{C}$ is closed under predecessors,
$add(\mathcal{C})$ is closed under submodules and, in particular,
under  kernels and images. Thus (1) and (2) are clearly
equivalent.

(2) is equivalent to  (3):  Since $add(\mathcal{C})$ is closed
under submodules, every quotient of a module in $add(\mathcal{C})$
is the  cokernel of a morphism in $add(\mathcal{C})$. Thus (2)
implies (3). The reverse implication is trivial.

 (3) implies (4): This is clear

(4) implies (5):  If $X\in\mathcal{C}$, then in  the radical
filtration $X\supset XJ(A)\supset XJ(A)^2\supset ....$ all the
terms are direct sums of predecessors of $X$. Hence, all  belong
to $add(\mathcal{C})$. Since $\mathcal{C}$ is closed under
predecessors, then, for every $k\geq 0$,   the projective cover
$P_k$ of $XJ(A)^k$ belongs to $add(\mathcal{C})$. The hypothesis 4
implies that $\frac{XJ(A)^k}{XJ(A)^{k+1}}\cong top(P_k)$ belongs
to $add(\mathcal{C})$. Since every composition factor of $X$ is
direct summand of some $\frac{XJ(A)^k}{XJ(A)^{k+1}}$, the
statement 5 follows.

(5) implies (2):   If $f:X\longrightarrow Y$ is a morphism between
modules in $\mathcal{C}$, the hypothesis guarantees that all
composition factors of $Z=coker(f)$ lie in $\mathcal{C}$. In
particular, $top(Z)\in add(\mathcal{C})$.  Since $\mathcal{C}$ is
closed under predecessors, we have $Z\in add(\mathcal{C})$, so
that $add(\mathcal{C})$ is closed under cokernels.

(3) implies (6)
:  Since  $\mathcal{C}$ is closed  under
predecessors, the projective cover of a module $X\in\mathcal{C}$
belongs to $add(\mathcal{C})$ and, consequently, to $add(P)$.
Hence $add(\mathcal{C})\subseteq Gen(P)$. The reverse inclusion
follows from the fact that
 $add(\mathcal{C})$ is closed under quotients.

Since (6) trivially implies (4), the proof is complete.
\end{proof}

We recall that an additive full subcategory $\mathcal{D}$ of
$mod_A$ is {\bf contravariantly finite} if, for every $X\in
mod_A$, there is a morphism $f:D_X\longrightarrow X$ (called a
right approximation) such that $D_X\in\mathcal{D}$ and, for any
other morphism $g:D\longrightarrow X$, with $D\in\mathcal{D}$,
there exists $h:D\longrightarrow D_X$ such that $f\circ h=g$. {\bf
Covariantly finite} subcategories are defined dually, and a
subcategory  is called {\bf functorially finite} if it is both
covariantly and contravariantly finite (see \cite{AS}).

\begin{cor} \label{dos-consecuencias}
Let $\mathcal{C}$ be a full subcategory of $ind_A$ closed under
predecessors such that $add(\mathcal{C})$ is an abelian exact
subcategory of $mod_A$, and let $e\in A$ be an idempotent such
that $P=eA$ is (isomorphic to) the supporting projective of
$\mathcal{C}$. The following assertions hold:

\begin{enumerate}
\item[(1)]$add(\mathcal{C})$ is functorially finite in $mod_A$
\item[(2)]$(Gen((1-e)A),add(\mathcal{C}))=(Gen((1-e)A),Gen(eA))$
is the split torsion pair in $mod_A$ having $add(\mathcal{C})$ as
torsion-free class. \item[(3)]The   torsion radical $t$ of the
above torsion pair is given by $t(X)=$ \newline $X(1-e)A$, for
every $X\in mod_A$
\end{enumerate}
\end{cor}
\begin{proof}

\begin{enumerate}
\item[(1)] Every torsion-free class is covariantly finite. By
Proposition \ref{auxilio-fundamental}(6),
$add(\mathcal{C})=Gen(P)$ is contravariantly finite, the (minimal)
right approximation of $X$ being the inclusion
$t_P(X)\hookrightarrow X$, where $t_P(X)$ is the trace of $P$ in
$X$, that is,  $t_P(X)=\sum_{f\in Hom_A(P,X)}Im(f)$. \item[(2)]
Consider the split torsion pair $(\mathcal{T},add(\mathcal{C}))$.
Since
 $\mathcal{C}$ is closed under composition factors, an $A$-module $X$
 lies in $\mathcal{T}$ if, and only if, $top(X)$ contains no simple
  summand from $\mathcal{C}$, that is, if and only if $top(X)\in Gen((1-e)A)$.
  This is equivalent to saying  that  $X$ is
  generated by $(1-e)A$. \item[(3)]  $t(X)$ is the
  (unique)
  maximal submodule of $X$  belonging to $Gen((1-e)A)$, which is the trace
    $t(X)=\sum_{f\in Hom_A((1-e)A,X)}Im(f)=X(1-e)A$
  \end{enumerate}
\end{proof}

We recall that an $A$-module $I$ is called {\bf hereditary
injective} if every  quotient of $I$ (or of $I^r$, with $r>0$)  is
an injective $A$-module.

\begin{rem} \label{trace}
 If $A=
\begin{pmatrix} C & 0\\ M & B \end{pmatrix}$, where $M$ is a
$B-C-$bimodule, then the right $A$-modules can be viewed
 as triples $(X,Y,\varphi )$, where $X\in mod_C$, $Y\in mod_B$ and
 $\varphi :Y\otimes_BM\longrightarrow X$ is a morphism in $mod_C$
 (see \cite{ARS}[Chapter III]). In this case, we may, and shall, identify
  $mod_C$
 with the full subcategory of $mod_A$ having as objects the triples
 $(X,0,0)$, with $X\in mod_C$.
\end{rem}

For any right artinian ring $R$, we  denote by  $gl.dim(R)$ the
global dimension of $R$. We are now able to state, and prove, the
main result of this paper.

\begin{teor}
\label{main-theorem} Let $A$ be a basic connected right artinian
ring and  $\mathcal{C}$ be a full subcategory of $ind_A$. The
following assertions are equivalent:

\begin{enumerate}
\item[(1)] $\mathcal{C}$ is closed under predecessors and
$add(\mathcal{C})$ is an abelian exact subcategory of $mod_A$
\item[(2)]There exists a ring isomorphism $A\cong\begin{pmatrix} C
& 0\\ M & B\\
\end{pmatrix}$ such that $M_C$ is a hereditary injective $C$-module and
$add(\mathcal{C})\cong mod_C$ \item[(3)] There exists an
idempotent $e\in A$ such that $eA(1-e)=0$, $\mathcal{C}$ consists
of those $X\in ind_A$ such that $Xe=X$ and every $Y\in
ind_A\setminus\mathcal{C}$ is generated by $(1-e)A$
\end{enumerate}

Further, if this is the case, then  $gl.dim(C)=Sup\{pd(X_A):$
$X\in\mathcal{C}\}$
\end{teor}
\begin{proof}

 (1) implies (3):  Let $e\in A$ be an idempotent such that $eA=P$ is
 the supporting projective of $\mathcal{C}$. By proposition
 \ref{auxilio-fundamental}, $add(\mathcal{C})=Gen(eA)$ and, by corollary
 \ref{dos-consecuencias},   the corresponding split torsion pair is
 $(Gen((1-e)A,Gen(eA))$. Therefore $eA(1-e)\cong
 Hom_A((1-e)A,eA)=0$ and so $XeA=Xe$, for all $X\in mod_A$. Hence
 $X\in add(\mathcal{C})$ if, and only if,  $X=Xe$. The last statement follows from the
 fact that the torsion pair is split.

\vspace*{0.5cm}
 (3) implies (1):  Since $eA(1-e)=0$, we have
 $Gen(eA)=\{X\in mod_A:$ $X=Xe\}$. The hypothesis (3) says
 exactly  that
 $(Gen((1-e)A),Gen(eA))=(Gen((1-e)A),add(\mathcal{C})$ is a split
 torsion pair. The statement then follow from lemma
 \ref{split-torsion} and proposition \ref{auxilio-fundamental}

\vspace*{0.5cm} (2) implies (3):  Setting $e=
\begin{pmatrix} 1 & 0\\ 0 & 0 \end{pmatrix} $, we have $1-e=
\begin{pmatrix} 0 & 0\\ 0 & 1 \end{pmatrix} $  so that,
clearly,
 $eA(1-e)=0$. The equality $\mathcal{C}=\{X\in ind_A:$
$Xe=X\}=\{X\in ind_A:$ $X(1-e)=0\}$ follows from the
interpretation of $mod_C$ as a full subcategory of $mod_A$.
  There remains to prove that
$ind_A\setminus\mathcal{C}\subset Gen(1-e)A)$. Let
$X\notin\mathcal{C}$ be indecomposable.  We claim that

\begin{center}
$X(1-e)AeA=X(1-e)A\cap XeA$ \end{center} Clearly, we have
$X(1-e)AeA\subseteq X(1-e)A\cap XeA$. Conversely, if $x\in
X(1-e)A\cap XeA$ then $x=xe$, due to the equality $XeA=Xe$. On the
other hand,
 $x=\sum_{1\leq i\leq n}y_i(1-e)a_i$, with $a_i\in A$ and $y_i\in
 X$. But then $x=xe=\sum_{1\leq i\leq n}y_i(1-e)a_ie$ belongs to
 $X(1-e)AeA$, thus establishing our claim.

 The $A$-module
$X(1-e)AeA$ is generated by $(1-e)AeA=M$ which, by hypothesis, is
a hereditary injective $C$-module. Hence, $X(1-e)AeA$ is injective
in $mod_C$,  and so we have a decomposition

\begin{center}
$XeA=X(1-e)AeA\oplus X'$ \end{center} in $mod_C$. Considering this
decomposition in $mod_A$ via the  embedding $mod_C\hookrightarrow
mod_A$, we have

\begin{center}
$X=XeA+X(1-e)A=X(1-e)AeA+X'+X(1-e)A=X'+X(1-e)A$ \end{center}  But
$X'\cap X(1-e)A\subseteq XeA\cap X(1-e)A=X(1-e)AeA$, and so
$X'\cap X(1-e)A\subseteq X'\cap X(1-e)AeA=0$. We thereby get a
decomposition

\begin{center}
$X=X'\oplus X(1-e)A$ \end{center} in $mod_A$. Since $X_A$ is
indecomposable and $X(1-e)\neq 0$ (because $X\notin\mathcal{C}$),
we conclude that $X'=0$ and, hence, $X=X(1-e)A\in Gen((1-e)A)$ as
desired.

\vspace*{0.5cm} (1) and (3) imply (2):  From (3), letting $C=eAe$,
$e'=1-e$,  $B=e'Ae'$ and $M=e'Ae$, we may identify $A$ with the
matrix algebra $A=
\begin{pmatrix} C & 0\\ M & B \end{pmatrix}$  and
$\mathcal{C}$ with $ind_C=\{X\in ind_A:$ $Xe=X\}=\{X\in ind_A:$
$Xe'=0\}$. By corollary \ref{dos-consecuencias}, the torsion
radical associated with the split torsion pair
$(Gen(e'A),add(\mathcal{C}))$ is given by $t(X)=Xe'A$, so that
$Xe'A$ is a direct summand of $X$, for every $X\in mod_A$. Let us
fix a complete set of primitive orthogonal idempotents
$\{e_1,...,e_s\}$ of $B=e'Ae'$, so that $e'=e_1+...+e_s$.
Interpreting $A$-modules as triples, as in  remark \ref{trace}(2),
we have that $e_iJ(A)e'A=(e_iJ(B)M,e_iJ(B),\mu )$ is a direct
summand of $e_iJ(A)=(e_iM,e_iJ(B),\mu ')$ in $mod_A$ (where $\mu$,
 $\mu'$ are the respective multiplication maps). This clearly
 implies that  $e_iJ(B)M$ is a direct
summand of $e_iM$ in $mod_C$.

We fix, for each $i$, a decomposition $e_iM=M'_i\oplus e_iJ(B)M$
in $mod_C$,  and let $M'=\oplus_{1\leq i\leq s}M'_i$. Then
$e_iM'=M'_i$ for each $i$, and $M=M'\oplus J(B)M$ in $mod_C$. An
easy induction shows that $J(B)M=J(B)M'+J(B)^kM$,
 for all
$k>0$. The nilpotency of  $J(B)$ yields $J(B)M=J(B)M'$. Since the
left multiplication by an element of $B$ gives an endomorphism of
$M_C$, the equality $J(B)M=J(B)M'$ implies that $J(B)M$ is
generated by $M'$ in $mod_C$. But then $M_C=M'\oplus J(B)M$  is
also generated by $M'$ in $mod_C$. In order to prove that $M_C$ is
hereditary injective, it suffices to show that each $M'_i=e_iM'_C$
is hereditary injective.

 Suppose that this is not the case and consider an
epimorphism $g:N\twoheadrightarrow Z$, where $N$ is an
indecomposable summand of some $M'_i$ and $Z$ is a non-injective
indecomposable $C$-module. Decomposing $M'_i=N\oplus N'$ in
$mod_C$, we see that $e_iJ(A)=N\oplus N'\oplus e_iJ(B)A$ is a
decomposition in $mod_A$, where the $C$-modules $N$, $N'$ are
viewed as $A$-modules and, by definition,
$e_iJ(B)A=(e_iJ(B)M,e_iJ(B),\mu )$, with $\mu$ the multiplication
map. We deduce an  embedding $Ker(g)\oplus N'\oplus
e_iJ(B)A\hookrightarrow e_iJ(A)\hookrightarrow e_iA$. The
corresponding quotient $X=\frac{e_iA}{Ker(g)\oplus N'\oplus
e_iJ(B)A}$ has simple top, hence is indecomposable. We also have

\begin{center}
$XJ(A)=\frac{e_iJ(A)}{Ker(g)\oplus N'\oplus
e_iJ(B)A}=\frac{N\oplus N'\oplus e_iJ(B)A}{Ker(g)\oplus N'\oplus
e_iJ(B)A}\cong Z$
\end{center}

Since $Z$ is not injective in $mod_C$, the functor $Ext_C^1(-,Z)$
is non-zero. It is easily seen that this is equivalent to the
existence of some simple $C$-module $S$ such that
$Ext_C^1(S,Z)\neq 0$. We fix a non-split exact sequence

\begin{center}
$0\rightarrow
Z\stackrel{j}{\longrightarrow}V\stackrel{p}{\longrightarrow}S\rightarrow
0$ \end{center}  in $mod_C$ which, clearly, is also non-split in
$mod_A$. By the above comments, the canonical inclusion
$XJ(A)\hookrightarrow X$ induces an embedding $i:Z\longrightarrow
X$. We thus have an amalgamated sum (pushout) diagram:

\vspace*{0.5cm}

\setlength{\unitlength}{1mm}
\begin{picture}(140,30)
\put(20,27){$0$} \put(22,28){\vector(1,0){5}} \put(28,27){$Z$}
\put(31,28){\vector(1,0){25}} \put(43,30){$j$} \put(58,27){$V$}
\put(61,28){\vector(1,0){25}} \put(73,30){$p$} \put(88,27){$S$}
\put(91,28){\vector(1,0){5}} \put(98,27){$0$}

\put(20,2){$0$} \put(22,3){\vector(1,0){5}} \put(28,2){$X$}
\put(31,3){\vector(1,0){25}} \put(43,5){$u$} \put(58,2){$W$}
\put(61,3){\vector(1,0){25}} \put(73,5){$w$} \put(88,2){$S$}
\put(91,3){\vector(1,0){5}} \put(98,2){$0$}

\put(29,26){\vector(0,-1){20}} \put(59,26){\vector(0,-1){20}}
\put(88,26){\line(0,-1){20}} \put(89,26){\line(0,-1){20}}
\put(27,16){$i$} \put(57,16){$r$}
\end{picture}

Since $(Gen(e'A),add(\mathcal{C}))$ is a split torsion pair, we
have $W=W_1\oplus W_2$, with $W_1\in add(\mathcal{C})=Gen(eA)$
(whence it is a $C$-module) and $W_2\in Gen(e'A)$. Since $X\in
Gen(e'A)$, the composition of $u$ with the projection
$W\longrightarrow W_1$ vanishes, so that $u(X)\subseteq W_2$. The
obvious inequalities between composition lengths $l(X)\leq
l(W_2)\leq l(W)=l(X)+1$ lead to two cases:

\begin{enumerate}
\item Assume first that $l(W_2)=l(W)=l(X)+1$. Then $W=W_2$ and
$W_1=0$, so that $W\in Gen(e'A)$. But $w:W\longrightarrow S$ is
non-zero, and $S\in\mathcal{C}$. This is a contradiction. \item
Assume $l(X)=l(W_2)=l(W)-1$. Identifying $X$ with $u(X)$, we have
$X=W_2$ so that $W_1\cong W/X\cong S$ and the above diagram
becomes

\vspace*{0.5cm}

\setlength{\unitlength}{1mm}
\begin{picture}(140,30)
\put(20,27){$0$} \put(22,28){\vector(1,0){5}} \put(28,27){$Z$}
\put(31,28){\vector(1,0){25}} \put(43,30){$j$} \put(58,27){$V$}
\put(61,28){\vector(1,0){25}} \put(73,30){$p$} \put(88,27){$S$}
\put(91,28){\vector(1,0){5}} \put(98,27){$0$}

\put(20,2){$0$} \put(22,3){\vector(1,0){5}} \put(28,2){$X$}
\put(31,3){\vector(1,0){21}} \put(40,5){$(1$ $0)^t$}
\put(54,2){$X\oplus S $} \put(65,3){\vector(1,0){21}}
\put(70,5){$(0$ $1)$} \put(88,2){$S$} \put(91,3){\vector(1,0){5}}
\put(98,2){$0$}

\put(29,26){\vector(0,-1){20}} \put(59,26){\vector(0,-1){20}}
\put(88,26){\line(0,-1){20}} \put(89,26){\line(0,-1){20}}
\put(27,16){$i$} \put(60,16){$(h$ $p)^t$}
\end{picture}

 for some  $h:V\longrightarrow X$. In particular, $h\circ
j=i$ and $p\circ j=0$. On the other hand, since $V$ is a
$C$-module, we have $Im(h)\subseteq XeA\subseteq XJ(A)\cong Z$
because $X$ has a simple top isomorphic to
$S_i=\frac{e_iA}{e_iJ(A)}$. We then get a morphism
$h':V\longrightarrow Z$ such that $i\circ h'=h$. But then
$i=h\circ j=i\circ h'\circ j$ and, since $i$ is a monomorphism, we
get $h'\circ j=1_Z$. This contradicts the fact that the upper
sequence in the above diagram is not split.
\end{enumerate}

In either case we have reached a contradiction. Hence each $M'_i$
is hereditary injective. That completes the proof of the
equivalence of (1), (2) and (3).

The last statement of the theorem follows from the fact that,  if
we identify $add(\mathcal{C})$ with the full subcategory $mod_C$
of $mod_A$, then the minimal projective resolution of any
$X\in\mathcal{C}$ is the same in $mod_C$ and $mod_A$.
\end{proof}

\vspace*{0.3cm}

Given a complete set of primitive orthogonal idempotents
$\mathcal{E}=\{e_1,...,e_n\}$ of $A$, and a subset $\Sigma
=\{e_{i_1},...,e_{i_r}\}$ of $\mathcal{E}$, we denote by
$e_\Sigma$ the sum $e_{i_1}+...+e_{i_r}$. With this notation, the
desired classification of the split torsion pairs with
torsion-free class closed under quotients follows directly from
our theorem.

\begin{cor} \label{clasification}
Let  $\mathcal{E}=\{e_1,...,e_n\}$ be a complete set of primitive
orthogonal idempotents of $A$. There is a one-to-one
correspondence between:

\begin{enumerate}
\item[(1)] The  full subcategories $\mathcal{C}$ of $ind_A$ closed
under predecessors such that $add(\mathcal{C})$ is an abelian
exact subcategory of $mod_A$ \item[(2)] The  split torsion pairs
in $mod_A$, with torsion-free class  closed under quotients
\item[(3)] The subsets $\Sigma\subseteq\mathcal{E}$ such that
$(1-e_\Sigma )Ae_\Sigma$ is a hereditary injective $e_\Sigma
Ae_\Sigma$-module and $e_\Sigma A(1-e_\Sigma )=0$
\end{enumerate}
\end{cor}

\section{Applications to Artin algebras}
Throughout this section, we assume that our algebras are basic and
connected  Artin algebras.
 We
denote by $Q_A$ the (valued) quiver of $A$ and by $(Q_A)_0$ the
set of points of $Q_A$. The idempotent corresponding to a point
$x\in (Q_A)_0$ is denoted by $e_x$, while we denote by $P_x$ (or
$S_x$) the corresponding indecomposable projective ( or simple,
respectively). For general facts about the module category of $A$,
we refer the reader to \cite{ARS}.

A first consequence of our main theorem is the following
combinatorial result:

\begin{cor} \label{B-C-connection}
Let $A$ be an algebra satisfying the equivalent conditions of the
theorem.  Then, for every arrow $y\rightarrow x$ in $Q_A$, with
$y\in (Q_B)_0$ and $x\in (Q_C)_0$, the point $x$ is a source in
$Q_C$.
\end{cor}
\begin{proof}
Since there exists an arrow $y\rightarrow x$ in $Q_A$, then
$\frac{e_yJ(A)e_x}{e_yJ(A)^2e_x}\neq 0$. Notice that  $e_yJ(A)e_x$
is identified with $e_yMe_x$ and $e_yJ(A)^2e_x$  with
$e_y[J(B)M+MJ(C)]e_x$. Then
$\frac{e_yMe_x}{e_y[J(B)M+MJ(C)]e_x}\neq 0$ and, in particular,
$\frac{e_yMe_x}{e_yMJ(C)e_x}\neq 0$. This  says that the simple
$C$-module $S_x$ is a direct summand of the top of the $C$-module
$e_yM$ and, hence, also of $top(M_C)$. Since $M_C$ is hereditary
injective, we conclude that $S_x$ is a simple injective
$C$-module, so that $x$ is a source in $Q_C$.
\end{proof}

We now consider the  case where $\mathcal{C}$ is the left part
$\mathcal{L}_A$ of $mod_A$, that is,  the full subcategory of
$ind_A$ consisting of those $X\in ind_A$ such that every
predecessor of $X$ has projective dimension at most one (see
\cite{HRS}). Thus, $\mathcal{L}_A$ is closed under predecessors.
The endomorphism algebra of the supporting projective of
$\mathcal{L}_A$ is denoted by $A_\lambda$ and is called the {\bf
left support} of $A$ (see \cite{ACT} and \cite{Sk}).

 We recall that $A$ is called {\bf left supported}
when $add(\mathcal{L}_A)$ is contravariantly finite in $mod_A$
(see \cite{ACT}). Many important classes of algebras are left
supported such as, for instance, the laura algebras which are not
quasi-tilted (see \cite{ACT}, \cite{Sk}).

\begin{cor} \label{LA}
Let $A$ be an Artin algebra such that $add(\mathcal{L}_A)$ is an
abelian exact subcategory of $mod_A$. Then:

\begin{enumerate}
\item[(1)] The left support $A_\lambda$ of $A$ is hereditary
\item[(2)] The algebra $A$ is left supported \item[(3)] If,
furthermore, the valued quiver of $A$ has no oriented cycles, then
$A=A_\lambda$. In particular, $A$  itself is  hereditary
\end{enumerate}
\end{cor}
\begin{proof}
(1) follows from the last statement of the theorem,  and (2)
follows from corollary \ref{dos-consecuencias}(1). In order to
prove (3), suppose that $A\neq A_\lambda$. There exists a point
$x_0\in (Q_A)_0$ such that $P_{x_0}\notin\mathcal{L}_A$. In
particular, the radical $P_{x_0}J(A)$ of $P_{x_0}$ admits an
indecomposable summand $R_{x_0}$ which is not in $\mathcal{L}_A$.
Hence there exists a point $x_1\in (Q_A)_0$ such that
$P_{x_1}\notin\mathcal{L}_A$ and $Hom_A(P_{x_1},R_x)\neq 0$. This
yields a non-zero non-isomorphism $f_1:P_{x_1}\longrightarrow
P_{x_0}$. Repeating the process for $x_1$ instead of $x_0$ yields
a point $x_2\in (Q_A)_0$ such that $P_{x_2}\notin\mathcal{L}_A$
and there exists a non-zero non-isomorphism
$f_2:P_{x_2}\longrightarrow P_{x_1}$. Inductively, we get a
sequence of non-zero non-isomorphisms between indecomposable
projective modules $...P_{x_n}\stackrel{f_n}{\longrightarrow}
P_{x_{n-1}}...\stackrel{f_2}{\longrightarrow}P_{x_1}
\stackrel{f_1}{\longrightarrow}P_{x_0}$. Since $(Q_A)_0$ is
finite, this sequence yields necessarily an oriented cycle in
$Q_A$, which is a contradiction.
\end{proof}

We note that, if $A_\lambda$ is hereditary, it does not follow in
general that $add(\mathcal{L}_A)$ is an abelian exact subcategory
of $mod_A$, as is shown by the following example.

\begin{ejem} \label{ejemplos1}
 Let $K$ be a field and  $A$
be the radical square zero $K$-algebra given by the quiver

\vspace*{0.5cm}

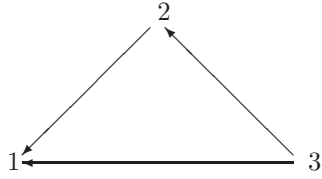
\setlength{\unitlength}{1mm}
\begin{picture}(140,20)
\put(30,0){$1$} \put(70,0){$3$} \put(68,1){\vector(-1,0){36}}
\put(50,20){$2$} \put(49,19){\vector(-1,-1){17}}
\put(68,2){\vector(-1,1){17}}
\end{picture}

Here $\mathcal{L}_A=\{P_1,P_2\}$ and its support is the hereditary
$K$-algebra with quiver
\begin{center}
$1\longleftarrow 2$ \end{center} However, $add(\mathcal{L}_A)$ is
not an abelian exact subcategory of $mod_A$ because it does not
contain the cokernel $S_2$ of the inclusion $P_1\longrightarrow
P_2$

\end{ejem}

 Our final application  is to local extensions of
 hereditary algebras. We recall that a triangular matrix algebra $A=
\begin{pmatrix} H & 0\\ M & R \end{pmatrix}$, where
$_RM_H$ is an $R-H-$bimodule, is called a {\bf local extension} of
$H$ in case $R$ is a local algebra (see \cite{MP}). Taking $R$ a
skew field, we see that this notion generalizes that of a
one-point extension. However, we are interested in the case where
$R$ is not a skew field, a hypothesis that we assume in the
sequel. We denote by $y$ the unique point in $Q_R$. For general
facts about the module category of a local extension, we refer the
reader to \cite{MP}

\begin{lema} \label{loc.ext.support}
Let $A= \begin{pmatrix} H & 0\\ M & R \end{pmatrix}$  be a local
extension of the hereditary algebra $H$. Then the left support
$A_\lambda$ is equal to $H$.
\end{lema}
\begin{proof}
Let $P_x$ be any indecomposable projective $H$-modules. The
predecessors of $P_x$ in $ind_A$  are (projective) $H$-modules
and, hence, $P_x\in\mathcal{L}_A$. On the other hand, the only
other indecomposable projective $P_y$ lies on an oriented cycle of
projectives in $ind_A$. Therefore $y\notin (Q_{A_\lambda})_0$,
because $A_\lambda$ is quasi-tilted by \cite{ACT}[2.1] and hence
triangular by \cite{HRS}
\end{proof}

It follows from the above lemma,  or from \cite{ACT}[2.1], that we
have an inclusion $\mathcal{L}_A\subseteq ind_H$. Our final result
says exactly when equality holds:

\begin{prop} \label{final}
Let $A=
\begin{pmatrix} H & 0\\ M & R \end{pmatrix}$  be a local
extension of the hereditary algebra $H$, where $R$ is not a
skew-field. The following statements are equivalent:

\begin{enumerate}
 \item[(1)] $add(\mathcal{L}_A)$ is an abelian exact subcategory of
$mod_A$ \item[(2)] $\mathcal{L_A}=ind_H$ \item[(3)] $M_H$ is
injective
\end{enumerate}

\end{prop}
\begin{proof}
(1) implies (3): By lemma \ref{loc.ext.support}, we have
$H=A_\lambda$. Our main theorem \ref{main-theorem} gives that
$M_H$ is injective.

(3) implies (2): From  theorem \ref{main-theorem} we get
$\mathcal{C}=ind_H$. Also, for any $X\in ind_H$, we have
$pd_A(X)\leq gl.dim (H)=1$. Then $ind_H\subseteq
add(\mathcal{L}_A)$, so that $\mathcal{L}_A=ind_H$

(2) implies (1):  The hypothesis  gives
$mod_H=add(\mathcal{L_A})$, and the statement follows at once.
\end{proof}

\begin{ejem}
Let $K$ be a field and let $A$ be the $K$-algebra given by the
quiver

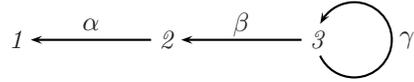
\begin{center}
\psset{unit=1mm}
\begin{pspicture}(0,0)(70,20)
\rput(0,10){1} \rput(30,12){$\beta$}
 \rput(40,10){3} \rput(20,10){2} \rput(10,12){$\alpha$}
 \rput(52,10){$\gamma$}
 \psline{->}(18,10)(2,10) \psline{->}(38,10)(22,10)
 \psarc[arrows=->](45,10){5}{-155}{155}
\end{pspicture}
\end{center}

 with relations $\gamma^2=0$
and $\gamma\beta\alpha =0$. Denoting the indecomposables by their
Loewy series, the regular module $A_A$ is given by:

\vspace*{0.5cm} \setlength{\unitlength}{1mm}
\begin{picture}(140,20)
\put(35,10){$1$} \put(40,10){$\oplus$} \put(45,13){$2$}
\put(45,7){$1$} \put(50,10){$\oplus$} \put(55,10){$2$}
\put(55,4){$1$} \put(58,16){$3$} \put(61,10){$3$} \put(61,4){$2$}
\end{picture}

Here,  $A$ is a local extension of the hereditary algebra $H$
given by the quiver $1\stackrel{\alpha}{\longleftarrow}2$, taking
$M_H=P_2\oplus S_2$, which is an injective $H$-module. The
hypothesis of proposition \ref{final} is satisfied, and therefore
$add(\mathcal{L}_A)=mod_H$ is an abelian exact subcategory of
$mod_A$. Notice that if we put here $C=H$ and
$B=R=K[\gamma]/(\gamma^2)$, then $J(B)M\cong S_2$ so that, taking
$M'=P_2$, we get the decomposition $M_C=M'\oplus J(B)M$ of the
proof of theorem \ref{main-theorem}.
\end{ejem}


\begin{thebibliography}{9999}

\bibitem{AF}{ANDERSON, F.W.; FULLER, K.R.}: "Rings and categories
of modules", 2nd edition. Springer-Verlag (1992).
\bibitem{ACT}{\sc ASSEM, I.; COELHO, F.U.; TREPODE, S.}: The left
and right parts of a module category. Preprint.
\bibitem{ARS}{\sc AUSLANDER, M.; REITEN, I.; SMALO, S.O.}:
"Representation theory of Artin algebras". Cambridge Univ. Press
{\bf 36}, in Cambridge Studies in Advanced Mathematics. Cambridge
(1995).
\bibitem{AS}{\sc AUSLANDER, M.; SMALO, S.O.}: Preprojective
modules over Artin algebras. J. Algebra {\bf 66}(1) (1980),
61-122.

\bibitem{HRS}{\sc HAPPEL, D.; REITEN, I.; SMALO, S.O.}: "Tilting in abelian
categories and quasitilted algebras". Mem. Amer. Math. Soc. {\bf
120}(575) (1996), 473-526.
\bibitem{MP}{\sc MARTINS, M.I.; DE LA PEÑA, J.A.}: On local
extensions of algebras. Comm. Algebra {\bf 27}(3) (1999),
1017-1031.
\bibitem{Sk}{\sc SKOWRONSKI, A.}: On Artin algebras with almost
all indecomposable modules of projective or injective dimension at
most one. Cent. Eur. J. Math. {\bf 1}(1) (2003), 108-122.

\end{thebibliography}
\end{document}